\documentclass[EJP]{ejpecpmod} 

\usepackage{color}
\usepackage{tikz}
\usepackage{graphicx}
\usepackage{subfigure}
\usepackage{wrapfig}
\usepackage{hyperref}

\SHORTTITLE{Uniform Avoidance Coupling} 

\TITLE{Uniform Avoidance Coupling of Simple Random Walks}

\AUTHORS{%
  EJ Infeld\footnote{Ryerson University, Canada.
    \EMAIL{evainfeld@ryerson.ca}}}


\AMSSUBJ{60J10} 

\SUBMITTED{} 
\ACCEPTED{} 

\VOLUME{0}
\YEAR{2012}
\PAPERNUM{0}
\DOI{vVOL-PID}

\ABSTRACT{We start by introducing \emph{avoidance coupling} of Markov chains, with an overview of existing results. We then introduce and motivate a new notion, \emph{uniform avoidance coupling}. We show that the only Markovian avoidance coupling on a cycle is of this type, and that uniform avoidance coupling of simple random walks is impossible on trees, and prove that it is possible on several classes of graphs. We also derive a condition on the vertex neighborhoods in a graph equivalent to that graph admitting a uniform avoidance coupling of simple random walks, and an algorithm that tests this with run time polynomial in the number of vertices. }

\begin{document}

\section{What is avoidance coupling?}

Avoidance coupling of Markov processes was first introduced in 2013 by Angel et al. \cite{AC} The concept is fairly simple: given several Markov chains on the same state space, can we somehow implement them together in such a way that no two of them ever occupy the same state, but if we only see one it is faithful to its original probability matrix, i.e. looks like the Markov chain we started with.

The details of this are easiest to explain on the example of two simple random walks. Suppose that we have two tokens on the same graph, and they move one after the other. Can we set up this process in such a way that they never collide, and yet if we can only see one of them its behavior will be indistinguishable from a simple random walk? 

The idea of \textquotedblleft indistinguishable from a simple random walk" is rather subtle. Clearly, from any state we can only go to a neighbor, but other than that finite time behavior is meaningless. Every finite sequence of allowed moves is bound to happen at some point in a simple random walk. So we need to consider infinite time behavior.

The observed process needs to not only have the correct (perceived) stationary distributions for each state, but also not have any (perceived) correlation between the history of the token and where it is going to go next. But if the coupling itself is not a Markov process, these might be hard to define. We often find that we run out of tools.

We will therefore focus on couplings that are themselves Markov processes. This allows us to easily define stationary distributions for both the coupling and the individual tokens. We will now look at a few definitions, before continuing with examples.\pagebreak

\begin{definition}
A \emph{coupling} of Markov chains is an implementation of the chains on a common probability space, in such a way that each chain, viewed separately, is faithful to its transition matrix. 
\end{definition}

\begin{definition}
Let $X_t$ and $Y_t$ be the state of Markov chains $X$ and $Y$ at time $t$, respectively. A coupling of $X$ and $Y$ is an \emph{avoidance coupling} if for all $t\in\mathbb{N}$, $X_t\not=Y_t$ and $X_{t+1}\not=Y_t$.
\end{definition}

Notice that Definition 1.2 makes it explicit that the states move in turn. Moving from state $t$ to $t+1$, token $X$ is not allowed to step on $Y_t$. But by the time token $Y$ moves, token $X$ is already on state $X_{t+1}$ and so that's the state that $Y$ is not allowed to step on.

Since we will talk about Markovian couplings, we are free to use tools such as stationary distributions and transition probabilities of the two states together. We define a Markovian avoidance coupling as follows.

\begin{definition}
A coupling is \emph{Markovian} if the probability distribution of $(X_{t+1},Y_{t+1})$ given history $(X_0,Y_0),\dots,(X_{t-1},Y_{t-1}),(X_t,Y_t)$ is identical to the probability distribution of $(X_{t+1},Y_{t+1})$ given history $(X_0',Y_0'),\dots,(X_{t-1}',Y_{t-1}'),(X_t,Y_t)$. In other words, the probability distribution of $(X_{t+1},Y_{t+1})$ depends only on the previous state $(X_t,Y_t)$.\footnote{This is different from a definition provided by Angel et al., where the distribution of $X_{t+1}$ depends on $(X_t,Y_t)$ and the distribution of $Y_{t+1}$ depends on $(X_{t+1},Y_t)$. We will discuss the relationship between the two definitions and the results achieved so far in the next section.} 

A \emph{Markovian avoidance coupling} on a state space $M$ is a Markov chain on $M\times M$ such that for any $x,y,z\in M$ the stationary probability of $(x,x)$ is 0, and if the stationary probability $(x,y)$ is positive, the transition probability of $(x,y)\rightarrow(y,z)$ is 0.\medskip

\noindent This can be extended to $k$ walks as follows. Let the $k$ simple random walks be $X^i_t$, and the positions of all tokens on the graph at time $t$ be $(X^1_t,\dots,X^k_t)$. Then the coupling is \emph{Markovian} if $(X^1_{t+1},\dots,X^k_{t+1})$ depends only on $(X^1_t,\dots,X^k_t)$.
\end{definition}

\begin{example} Consider an avoidance coupling of two simple random walks on a cycle $C_n$, $n\geq  4$, constructed as follows. Start the chain in a state $(x,y)$, such that $x\not= y$, $x\not\sim y$ (i.e. $x$ and $y$ are not adjacent.) Then flip a fair coin, the same coin for both tokens. If it's heads, move them both clockwise. If it's tails, move them both anticlockwise. They will never collide, but if we can only see one of them, it is doing a faithful simple random walk.\end{example}

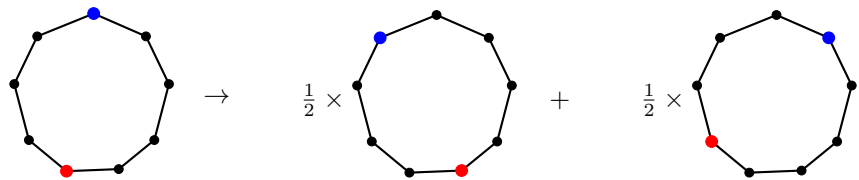
\begin{figure}[h!]\centering
\hfill\hfill
\begin{tikzpicture}[scale=0.55]
\foreach \x in {(0.05,2.05), (1.5,-1), (-1.5,-1), (1.3,1.5), (1.85,0.35), (0.65,-1.7), (-0.6,-1.75), (-1.3,1.5), (-1.85,0.35)}{
\filldraw \x circle (3pt);}
\draw[thick] (0.05,2.05) -- (1.3,1.5);
\draw[thick] (1.3,1.5) -- (1.85,0.35);
\draw[thick] (1.85,0.35)--(1.5,-1);
\draw[thick] (1.5,-1) -- (0.65,-1.7);
\draw[thick] (0.65,-1.7) -- (-0.6,-1.75);
\draw[thick] (-0.6,-1.75) -- (-1.5,-1);
\draw[thick] (-1.5,-1) -- (-1.85,0.35);
\draw[thick] (-1.85,0.35) -- (-1.3,1.5);
\draw[thick] (-1.3,1.5) -- (0.05,2.05);
\filldraw[blue] (0.05,2.05) circle (4pt);
\filldraw[red] (-0.6,-1.75) circle (4pt);
\draw (3,0) node {$\rightarrow$};
\end{tikzpicture}\hfill
\begin{tikzpicture}[scale=0.55]
\foreach \x in {(0.05,2.05), (1.5,-1), (-1.5,-1), (1.3,1.5), (1.85,0.35), (0.65,-1.7), (-0.6,-1.75), (-1.3,1.5), (-1.85,0.35)}{
\filldraw \x circle (3pt);}
\draw[thick] (0.05,2.05) -- (1.3,1.5);
\draw[thick] (1.3,1.5) -- (1.85,0.35);
\draw[thick] (1.85,0.35)--(1.5,-1);
\draw[thick] (1.5,-1) -- (0.65,-1.7);
\draw[thick] (0.65,-1.7) -- (-0.6,-1.75);
\draw[thick] (-0.6,-1.75) -- (-1.5,-1);
\draw[thick] (-1.5,-1) -- (-1.85,0.35);
\draw[thick] (-1.85,0.35) -- (-1.3,1.5);
\draw[thick] (-1.3,1.5) -- (0.05,2.05);
\filldraw[blue] (-1.3,1.5) circle (4pt);
\filldraw[red] (0.65,-1.7) circle (4pt);
\draw (3,0) node {$+$};
\draw (-3,0) node {$\frac{1}{2}$};
\draw (-2.4,0) node {$\times$};
\end{tikzpicture}\hfill\begin{tikzpicture}[scale=0.55]
\foreach \x in {(0.05,2.05), (1.5,-1), (-1.5,-1), (1.3,1.5), (1.85,0.35), (0.65,-1.7), (-0.6,-1.75), (-1.3,1.5), (-1.85,0.35)}{
\filldraw \x circle (3pt);}
\draw[thick] (0.05,2.05) -- (1.3,1.5);
\draw[thick] (1.3,1.5) -- (1.85,0.35);
\draw[thick] (1.85,0.35)--(1.5,-1);
\draw[thick] (1.5,-1) -- (0.65,-1.7);
\draw[thick] (0.65,-1.7) -- (-0.6,-1.75);
\draw[thick] (-0.6,-1.75) -- (-1.5,-1);
\draw[thick] (-1.5,-1) -- (-1.85,0.35);
\draw[thick] (-1.85,0.35) -- (-1.3,1.5);
\draw[thick] (-1.3,1.5) -- (0.05,2.05);
\filldraw[blue] (1.3,1.5) circle (4pt);
\filldraw[red] (-1.5,-1) circle (4pt);
\draw (-3,0) node {$\frac{1}{2}$};
\draw (-2.4,0) node {$\times$};
\end{tikzpicture}\hfill\hfill\hfill
\caption{Transition probabilities for a state in an avoidance coupling of two simple random walks on $C_9$.}
\end{figure}

It is not hard to see that for any $C_n$ we can construct an avoidance coupling of $k$ simple random walks for any $k\leq n/2$. Just start the tokens in states of which no two are adjacent and follow this algorithm!

\begin{example}Let's look at an example of a graph that does not admit an avoidance coupling of two simple random walks. Consider a path $P_n$, for any integer $n$, with states labelled from $0$ to $n-1$, from left to right. Suppose for contradiction, that such coupling is possible. But consider any state of such coupling, say $(i,j)$ (i.e. when one token is at $i$ and the other is at $j$) with $i<j$. Then there is no way the token that is currently at $i$ will ever reach $n-1$ without colliding with the other token. 

\vfill

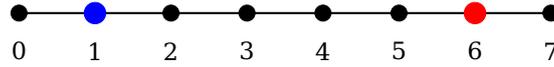
\begin{figure}[h!]\centering
\begin{tikzpicture}[scale=1]
\foreach \x in {0,1,2,3,4,5,6,7}{
\filldraw (\x,0) circle (3pt);}
\draw[thick] (0,0) -- (7,0);
\filldraw[blue] (1,0) circle (4pt);
\filldraw[red] (6,0) circle (4pt);
\foreach \x in {0,1,2,3,4,5,6,7}{
\draw (\x,-0.5) node {\x};}
\end{tikzpicture}
\caption{The blue token will not get to the vertex 7 without colliding with the red token.}
\end{figure}
The longer we watch this process the more we will be convinced that it is not in fact a simple random walk, since we expect a simple random walk to reach each state of the path.\end{example}

This is where the fact that the two tokens move in turn, rather than simultaneously, is critical. If the tokens were allowed to swap states, all bipartite graphs, including a path, would admit an avoidance coupling of simple random walks - just start them out on opposite sides of a bipartition and let each do a simple random walk! Later, when we talk about \emph{uniform avoidance coupling}, a stronger notion, we will be able to talk about the tokens as if they were moving simultaneously and it will still imply that a Markovian avoidance coupling of simple random walks exists.  We will also show that we can find such Markovian avoidance coupling on any bipartite graph that has no vertex of degree 1.

\section{Notation}

Unless stated otherwise, we will look at couplings of simple random walks on finite, connected graphs. For a graph $G=(V,E)$, let the state space be $V^k=V\times\dots\times V$, where  $V$ is the vertex set of $G$. For $v,w\in V$, let $v\sim w$ mean that $v$ and $w$ are adjacent, and $v\not\sim w$ mean that they are not.

 Consider a Markovian process $\mathcal{M}$ on state space $V\times\dots\times V$. Let $(v^1_{t},\dots,v^k_{t})$, $v^i_t\in V$ for all $1\leq i\leq k$ denote a state of this process, interpreted as token $1$ being at state $v^1_t$ and token $k$ being at state $v^k_t$ at time $t$.\medskip

\noindent A transition in this process will correspond to all tokens  making one move each: $$\mathcal{M}_t\rightarrow \mathcal{M}_{t+1}$$ $$(v^1_{t},\dots,v^k_{t})\rightarrow(v^1_{t+1},\dots,v^k_{t+1}).$$ For $v^i_t,v^i_{t+1}\in M$ for all $1\leq i\leq k$. 

\noindent Define the following:
\begin{itemize}
\item Let $s(v^1_{t},\dots,v^k_{t})$ be the stationary probability of the state $(v^1_{t},\dots,v^k_{t})$. \item Let $T[((v^1_{t},\dots,v^k_{t})\rightarrow(v^1_{t+1},\dots,v^k_{t+1})]$ be the transition probability from $(v^1_{t},\dots,v^k_{t})$ to $(v^1_{t+1},\dots,v^k_{t+1})$.\end{itemize}

\noindent This process is an \emph{avoidance} process if and only if:\begin{itemize}
\item $v^i_t=v^j_t$ for any $i\not=j$ $\Rightarrow$ $s(v^1_{t},\dots,v^k_{t})=0$
\item $s(v^1_{t},\dots,v^k_{t})\not=0$ and $v^i_t=v^j_{t+1}$ for some $i<j$ $\ \Rightarrow T[((v^1_{t},\dots,v^k_{t})\rightarrow(v^1_{t+1},\dots,v^k_{t+1})]=0$. In other words, no token that moves earlier can step on the previous position of a token that moves later.
\end{itemize}

\noindent If we only look at two tokens, these conditions can be written as:
\begin{itemize}
\item $s(v,v)=0,\ \forall v\in V$
\item $s(v,w)\not=0\ \Rightarrow T[(v,w)\rightarrow(w,u)]=0$, $\forall v,w,u\in V$. 

Note that we may have $T[(v,w)\rightarrow (u,v)]\not=0.$
\end{itemize}

Making sure that an avoidance Markov process is an avoidance coupling of simple random walks is much more complicated. For example, take a process with two tokens $X$ and $Y$ on $V\times V$. \medskip

\noindent Let $s(v),\ v\in V$ be the stationary probability of a simple random walk at vertex $v$. Then we need: $$\sum_{w\in V}s(v,w)=s(v).$$ But we also need to make sure that the transition probabilities for each token are independent of that token's history and correspond to the ones in a simple random walk. 

So, say, given that the token $X$ is at vertex $v$ and given any valid history of $X$, the probability of $v$ stepping to each adjacent vertex has to be the same. We will use this property in section 6 to show that the coupling in example 1.4 is in fact the only Markovian avoidance coupling of simple random walks on a cycle. 

In Example 2.1 we construct a Markov process of two tokens on a tree, in which the position of each token has the stationary distribution of a simple random walk, but the transition probabilities are not independent of history.

\begin{example}[Not a coupling] Consider the following tree, where each vertex is labelled with $(x,y)$ where $x$ is the label of the branch and $y$ is the vertex's depth:

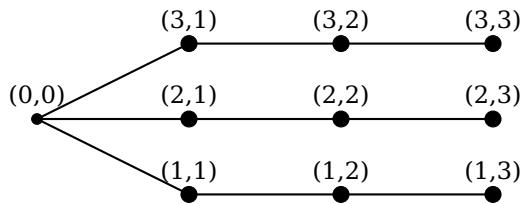
\begin{figure}[ht!]\centering
\begin{tikzpicture}[scale=1]
\filldraw (-2,1) circle (2pt);
\foreach \x in {0,1,2}{
\draw[thick] (-2,1) -- (0,\x) -- (4,\x);
\foreach \y in {0,2,4}{
\filldraw (\y,\x) circle (3pt);
}
}
\foreach \x in {1,2,3}{
\foreach \y in {1,2,3}{
\draw (2*\y-2,\x-0.7) node {(\x,\y)};
}
}
\draw (-2,1.3) node {(0,0)};
\end{tikzpicture}
\end{figure}
The tokens move according to following rules: whenever one token is at the root, the other is at a leaf and the sum of depths of the two tokens is always 3. The transition probabilities are:
\begin{itemize}
\item If one of the tokens is at $(0,0)$ (\emph{center}) and the other is at $(i,3)$, then the token at $(0,0)$ goes to each of $(j,1)$ with probability 1/3 and the token at $(i,3)$ goes to $(i,2)$.
\item If the two tokens are at $(i,1)$ and $(i,2)$ on some branch $i$, they move apart, i.e. go to $(0,0)$ and $(i,3)$.
\item If the two tokens are at $(i,1)$ and $(j,2)$ where $i\not=j$, they go to $(i,2)$ and $(j,1)$ with probability $3/4$ and to $(0,0)$ and $(j,3)$ with probability 1/4.
\end{itemize}
Each token in this process has the correct stationary distribution of a simple random walk:
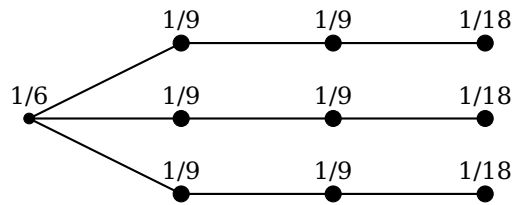
\begin{figure}[ht!]\centering
\begin{tikzpicture}[scale=1]
\filldraw (-2,1) circle (2pt);
\foreach \x in {0,1,2}{
\draw[thick] (-2,1) -- (0,\x) -- (4,\x);
\foreach \y in {0,2,4}{
\filldraw (\y,\x) circle (3pt);
}
}
\foreach \x in {1,2,3}{
\foreach \y in {1,2}{
\draw (2*\y-2,\x-0.7) node {1/9};
}
}
\foreach \x in {0,1,2}{
\draw (4,\x+0.3) node {1/18};
}
\draw (-2,1.3) node {1/6};
\end{tikzpicture}
\end{figure}

But the transition probabilities are not independent of history. Over all, if we see a token at $(1,1)$, it will go to $(1,2)$ with probability 1/2. But if we see it go from $(1,3)$ to $(1,2)$ and then $(1,1)$ then we can be sure that the other token is not on the same branch, and so it will go back to $(1,2)$ with probability $3/4$.\end{example}

\section{The literature so far: avoidance coupling on complete graphs.}

There are, as of writing this thesis, two published articles on this subject: the original paper by Angel et al.\ \cite{AC}, and a 2015 follow-up by Feldheim \cite{MAC}. They deal with avoidance coupling on the complete simple graph $K_n$, and the complete graph with loops $K_n^*$. We will now quickly discuss those of their results that are relevant to this paper.

\begin{figure}[ht!]\centering
\hfill
\begin{tikzpicture}[scale=2]
\path[white] (0,0.85) edge [loop] node {} (0,0.85);
\path[white] (-0.5,0) edge [in=185, out=275, loop] node {} (-0.5,0);
\path[white] (0.5,0) edge [in=265, out=355, loop] node {} (0.5,0);
\filldraw (0,0.85) circle (1pt);
\filldraw (-0.5,0) circle (1pt);
\filldraw (0.5,0) circle (1pt);
\draw (-0.5,0) -- (0.5,0) -- (0,0.85) -- (-0.5,0);
\end{tikzpicture}\hfill
\begin{tikzpicture}[scale=2]
\filldraw (0,0.85) circle (1pt);
\filldraw (-0.5,0) circle (1pt);
\filldraw (0.5,0) circle (1pt);
\draw (-0.5,0) -- (0.5,0) -- (0,0.85) -- (-0.5,0);
\path (0,0.85) edge [loop] node {} (0,0.85);
\path (-0.5,0) edge [in=185, out=275, loop] node {} (-0.5,0);
\path (0.5,0) edge [in=265, out=355, loop] node {} (0.5,0);
\end{tikzpicture}\hfill\hfill
\caption{Graphs $K_3$ and $K_3^*$, respectively.}
\end{figure}
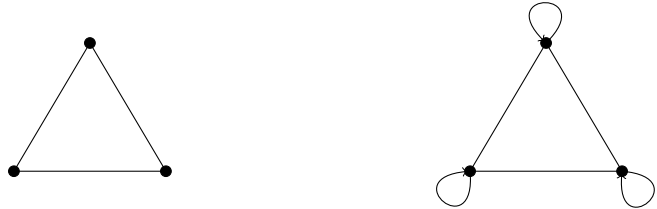

Both of these papers define Markovian coupling differently than we did in Definition 1.3. Angel et al.\ \cite{AC} use the strongest definition, and we will call those walks \emph{super-Markovian}. They introduce the concept of avoidance coupling and show that for any $n$, there exists an avoidance coupling of $k$ simple random walks on $K_n^*$ for any $k\leq n/4$, and a super-Markovian avoidance coupling on $K_n$ for any $k\leq n/(56\log_2n)$ and on $K_n^*$ for any $k\leq n/(8\log_2n)$.

In particular, we can construct a super-Markovian avoidance coupling of $k$ simple random walks on $K_n^*$, where $k\leq 2^d$ and $n=2^{d+1}$, and \cite{AC} proves monotonicity for $K^*_n$, i.e. shows that if there exists an avoidance coupling of $k$ walks on $K_n^*$, then there there is an avoidance coupling of $k$ walks on $K^*_{n+1}$. However, the monotonicity does not necessarily preserve the Markovian property. These two facts together imply that for any $n$ there exists an avoidance coupling of $k$ simple random walks on $K^*_n$ for any $k\leq n/4$.

Feldheim \cite{MAC} extends the monotonicity result to $K_n,$ and his construction preserves a weaker definition Markovian property. He calls these processes \emph{label-Markovian}.

\begin{definition}[Super-Markovian avoidance coupling of two simple random walks on $G$, \cite{AC}]Recall that on a graph $G=(V,E)$, an avoidance coupling of two simple random walks $X_{t,\ t\geq0}$ and $Y_{t,\ t\geq0}$ is an implementation of the two walks on a common probability space such that $\forall t\geq0$, $X_{t+1}\not=Y_t$ and $Y_{t+1}\not=X_{t+1}$. This means such that each, viewed separately, is faithful to its original probability matrix. This coupling is \emph{super-Markovian} if the probability distribution of $X_{t+1}$ is independent of history other than $X_t$ and $Y_t$, and the distribution of $Y_{t+1}$ is independent of history other than $X_{t+1}$ and $Y_t$. The state that the chain is in must contain a record of which token's move it is.\medskip

\noindent This definition can be extended to $k$ walks as follows. Let the $k$ simple random walks be $X_{t}^{i}$, with$,\ t\geq 0$ and $1\leq i\leq k$. An avoidance coupling is one where $\forall t\geq 0,$ $X_{t}^i\not=X_{t}^j$ for any $i\not= j$ and $X_{t+1}^i\not=X_{t}^j$ for any $i< j$. It's super-Markovian if $X^i_{t+1}$ depends only on $(X^1_{t+1},\dots,X^{i-1}_{t+1},X^i_t,\dots,X^k_t)$. The chain information must contain a record of which token's move it is.\end{definition}

In Definition 1.3, $(X_{t+1}^1,\dots,X_{t+1}^k)$ depends on $(X_t^1,\dots,X_t^k)$, and the moves of all tokens count as a single transition. This is clearly a weaker definition.

Any super-Markovian process is also Markovian. If we have $X_t$ and $Y_t$, an algorithm to generate $X_{t+1}$ out of those and an algorithm to generate $Y_{t+1}$ out of $X_{t+1}$ and $Y_{t}$, we can combine these algorithms to generate $(X_{t+1},\ Y_{t+1})$. Here is what this looks like more formally:

\begin{theorem}
Any super-Markovian coupling is also Markovian. That is, it has a form such that $(X_{t+1},Y_{t+1})$ depends only on $(X_t,Y_t)$, and the moves of both tokens count as a single transition.
\end{theorem}

\noindent\emph{Proof:} Let a super-Markovian chain be defined on some state space $V$, with two tokens $X$ and $Y$ moving in turn, such that $X_{t+1}$ depends only on $X_t$ and $Y_t$ and $Y_{t+1}$ depends only on $X_{t+1}$ and $Y_t$. Then for any states $X_t,\ Y_t, X_{t+1},\ Y_{t+1}$ such that if $X$ is at $X_t$ and it's $X$'s move, there is a non-zero probability that $Y$ is at $Y_t$, and if $Y$ is at $Y_t$ and it's $Y$'s move there is a non-zero probability that $X$ is at $X_{t+1}$. Define: \begin{itemize}\item  $f(X_t,Y_t,X_{t+1})$ to be the probability that while it's $X$'s move and the current positions are $X_t$ and $Y_t$, $X$ moves to $X_{t+1}$. \item  $g(X_{t+1},Y_t,Y_{t+1})$ to be the probability that while it's $Y$'s move and the current positions are $X_{t+1}$ and $Y_t$, $Y$ moves to $Y_{t+1}$.\end{itemize}

\noindent Define the transition probability as: $$T[(X_t,Y_t)\rightarrow(X_{t+1},Y_{t+1})]=f(X_t,Y_t,X_{t+1})g(X_{t+1},Y_t,Y_{t+1})$$
We can easily verify that: $$\sum_{Y_{t+1}}f(X_t,Y_t,X_{t+1})g(X_{t+1},Y_t,Y_{t+1})=f(X_t,Y_t,X_{t+1})
\sum_{Y_{t+1}}g(X_{t+1},Y_t,Y_{t+1})$$ $$=f(X_t,Y_t,X_{t+1})$$ And that regardless of $X_t$, given $X_{t+1}$ and $Y_t$, the move of $Y$ is governed by $g$. $\Box$
\medskip

\noindent We will now look at an example from \cite{AC} of a super-Markovian avoidance coupling on $K_n$ where $n$ is a composite integer. This process is not what we will later call a \emph{uniform} avoidance coupling.

\begin{example}[\cite{AC}]
Consider the complete graph $K_n$, where $n$ is a composite number, $n=ab,\ a,b>1$. There exists the following super-Markovian avoidance coupling of two simple random walks on this graph. 

First, divide the vertices of $K_n$ into $b$ clusters of size $a$. Start with tokens $X$ and $Y$ in two distinct clusters, and it's $X$'s move. 

\begin{itemize}
\item If it's $X$'s move, go to $Y$'s cluster with probability $\frac{a(b-1)}{ab-1}$. Choose an available vertex in that cluster uniformly at random. Otherwise, with probability $\frac{a-1}{ab-1}$ choose an available vertex in $X$'s own cluster uniformly at random.
\item If it's $Y$'s move, see if $X$ is in the same cluster. If yes, choose a vertex in a different cluster uniformly at random. If not, choose an available vertex in $Y$'s cluster uniformly at random.
\end{itemize}

To show that this is a coupling of simple random walks, it's helpful to look at the moves of both tokens together. $Y$ is moving with the probabilities of a simple random walk, since it effectively chooses any vertex other than the one it's at with uniform probability, and it's up to $X$ to move out of the way. Meanwhile, $X$ somewhat mimics the moves of $Y$. If $Y$ stays in the same cluster, so does $X$. If it doesn't, $X$ goes to $Y$'s previous cluster. Therefore, the sequence of clusters $X$ is in is as in a simple random walk. Within each cluster, $X$ can go to any vertex other than the one $Y$ is at, and $Y$ can be in at any of them with uniform probability. Between any two instances of $X$ changing clusters, any memory of which vertex within a cluster $Y$ was at is lost. Therefore $X$ is also doing a faithful simple random walk.

Using notation from Theorem 3.2, we can find the functions $f$ and $g$. Let $x\parallel y$ mean that $x$ and $y$ are distinct vertices in the same cluster, and $x{\nparallel} y$ mean that they are in distinct clusters. We have:

$$f(X_t,Y_t,X_{t+1})=\begin{cases}\frac{a(b-1)}{ab-1}\frac{1}{a-1} \text{\ if }Y_t\parallel X_{t+1} \\ \frac{a-1}{ab-1}\frac{1}{a-1} \text{\ \ \  if } X_t\parallel X_{t+1} \\ 0 \text{\ \ \ \ \ \ \ \ \ \ \ otherwise. } \end{cases}$$  $$g(X_{t+1},Y_t,Y_{t+1})=\begin{cases}\frac{1}{a(b-1)} \text{\ \ \ if }Y_t\parallel X_{t+1},\ Y_t\nparallel Y_{t+1} \\ \frac{1}{a-1} \text{\ \ \ \ \ \  if } Y_t\nparallel X_{t+1},\ Y_t\parallel Y_{t+1} \\ 0 \text{ \ \ \ \ \ \ \ \  otherwise. } \end{cases}$$ Where $X_t, X_{t+1}, Y_t$ are distinct vertices, and $Y_{t+1}$ is distinct from $X_{t+1}$ and $Y_t$, and where $X_t\nparallel Y_t$.\medskip

\noindent For all $X_t,\ Y_t$ such that $X_t\nparallel Y_t$, we get:
$$T[(X_t,Y_t)\rightarrow(X_{t+1},Y_{t+1})]=\begin{cases}\frac{1}{ab-1}\frac{1}{a-1}\frac{1}{a(b-1)} \text{\ \ \ if }Y_t\parallel X_{t+1},\ Y_t\nparallel Y_{t+1} \\ \frac{a-1}{ab-1}\frac{1}{a-1}\frac{1}{a-1} \text{\ \ \ \ \ \  if } X_t\parallel X_{t+1},\ Y_t\parallel Y_{t+1} \\ 0 \text{\ \ \ \ \ \ \ \ \ \ \ \ \ \ \ \ \ \ \  otherwise. } \end{cases}$$ In other words, 
$$T[(X_t,Y_t)\rightarrow(X_{t+1},Y_{t+1})] \begin{cases}\frac{1}{(ab-1)(a-1)} \text{\ \ \ if }(Y_t\parallel X_{t+1},\ Y_t\nparallel Y_{t+1})\text{ or }  (X_t\parallel X_{t+1},\ Y_t\parallel Y_{t+1}) \\ 0 \text{\ \ \ \ \ \ \ \ \ \ \ \ \ \  otherwise. } \end{cases}$$
\end{example}

\noindent The fact that some Markov chains admit a Markovian avoidance coupling but not a super-Markovian one can be illustrated by the following example from \cite{AC}.

\begin{theorem}[\cite{AC}] Let $M$ be a Markov chain with states $\{0,1,2\}$ and transition matrix: $$T=\begin{array}{c}
0 \\ 1 \\ 2
\end{array}\left[\begin{array}{ccc}\frac{1}{3} & \frac{1}{3} & \frac{1}{3} \\ \frac{1}{3} & \frac{1}{3} & \frac{1}{3} \\ \frac{1}{3} & \frac{1}{3} & \frac{1}{3}\end{array}\right].$$ Then a Markovian avoidance coupling of two such Markov chains is possible, but a super-Markovian avoidance coupling is not.\end{theorem}

\noindent\emph{Proof:} Consider the following Markovian avoidance process. Start the chain in any state $(X_0,\ Y_0)$ such that $X_0\not=Y_0$, with uniform probability. Suppose that the state at time $t$ is $(X_t,\ Y_t)\in M\times M,\ X_t\not= Y_t$. Then let $(X_{t+1},\ Y_{t+1})$ be a uniformly random pair of states such that $X_{t+1}\not=Y_t,\ Y_{t+1}\not= X_{t+1}$ and $(X_t,\ Y_t)\not=(X_{t+1},\ Y_{t+1})$.

\begin{wrapfigure}[7]{r}[-25pt]{2.7cm}
\raisebox{0pt}[\dimexpr\height-1\baselineskip\relax]{
\begin{tikzpicture}[scale=2]
\filldraw (0,0.85) circle (1pt);
\filldraw (-0.5,0) circle (1pt);
\filldraw (0.5,0) circle (1pt);
\draw (-0.5,0) -- (0.5,0) -- (0,0.85) -- (-0.5,0);
\path (0,0.85) edge [loop] node {} (0,0.85);
\path (-0.5,0) edge [in=185, out=275, loop] node {} (-0.5,0);
\path (0.5,0) edge [in=265, out=355, loop] node {} (0.5,0);
\draw (0,0.65) node {0};
\draw (-0.58,0.1) node {2};
\draw (0.58,0.1) node {1};
\end{tikzpicture}}
\end{wrapfigure}
Consider this process as two tokens walking on $K^*_3$ (pictured on the right.) We would like to show that this is a coupling of two simple random walks. 

For example, if $(X_t,Y_t)=(0,1)$ the next state is one of $\{(0,2),\ (2,0),\ (2,1) \}$. It is easy to see that at each point each token has the correct probability of staying in the same state, and that token $Y$ goes to each state with probability $\frac{1}{3}$ regardless of the current state $(X_t,Y_t)$. 

Is token $X$ performing a simple random walk as well? If while $X$ is at, say, $0$ the probability of $Y$ being at $1$ and $2$ is $\frac{1}{2}$ each, then $X$ is equally likely to step to either 1 or 2. We can check that this is true by induction. Suppose that at time $t$, given $X_t$, $Y_t$ can be either of the other two states with probability $\frac{1}{2}$. If $X_{t+1}=X_{t}$, $Y_t$ changes position. But the probability distribution stays the same. Otherwise, we observe $X$ changing positions. But given that knowledge, $Y_{t+1}$ is either where it was before or at $X_t$ with equal probability. So in both cases we have no extra knowledge about $Y_{t+1}$. We can conclude that this is a Markovian avoidance coupling of simple random walks.

However, a super-Markovian coupling does not exist. We follow the proof from \cite{AC}. Let $p_{x,y}$ be the probability that $X$ stays at $x$ given that it is $X$'s move and the current state is $(x,y)$. Let $q_{x,y}$ be the probability that $Y$ stays at $y$ given that the current position is $(x,y)$.

Suppose that $X$ has just moved from $0$ to $1$. The conditional probability of next moving back to $0$ is 1/3. $Y$ must have been at $2$ and will stay there with probability $q_{1,2}$ after which $X$ moves with probability $1-p_{1,2}$. So we must have $1/3=q_{1,2}(1-p_{1,2})$.

Similarly, suppose $Y$ just moved from $0$ to $2$. It will go back with probability $1/3=p_{1,2}(1-q_{1,2})$. So we have $p_{1,2}=q_{1,2}$, and similarly: $$p_{x,y}=q_{x,y},\ \forall x,y\in \{0,1,2\},\ x\not=y.$$ But the equation $p(1-p)=1/3$ has no real solutions. \hfill $\Box$\medskip

\noindent In general, \cite{AC} shows that for a Markov chain $M_s$, with the following transition matrix: $$T=\begin{array}{c}
0 \\ 1 \\ 2
\end{array}\left[\begin{array}{ccc}s & \frac{1-s}{2} & \frac{1-s}{2} \\ \frac{1-s}{2} & s & \frac{1-s}{2} \\ \frac{1-s}{2} & \frac{1-s}{2} & s\end{array}\right],\ s\in [0,1),$$ a Markovian avoidance coupling exists if and  only if $s\geq 1/3$ and a super-Markovian avoidance coupling exists if and only if $s\geq 1/2.$ 

For $s=0$, this is just a simple random walk on $K_3$. It is easy to see that no avoidance coupling of two of those is possible - at each point there is only one vertex to go to, and the tokens would either have to always go clockwise, or always go anti-clockwise.

We refer the reader to \cite{AC, MAC} for the remainder of existing results.

\section{Uniform avoidance coupling of simple random walks}

We said that the requirement that a token's behavior is indistinguishable from a simple random walk is a subtle one. It needs to not only have the correct stationary distribution at the states, but also transition probabilities must be independent of each token's history, other than the current state. The problem is, if we are coupling two processes, the behavior of one token might force a correlation between the other token's history and next move. To show a coupling is not an avoidance coupling of simple random walks, we can find such a correlation.

In order to find classes graphs that admit avoidance couplings of simple random walks, we introduce the stronger notion of \emph{uniform avoidance coupling} - a Markovian avoidance coupling that is faithful at every step. At each state, each token has the correct transition probabilities of going to every adjacent vertex, regardless of where the other token is. However, the steps of the two tokens are not independent. 

The avoidance coupling on a cycle in Example 1.4 is a uniform avoidance coupling of simple random walks.

\begin{definition}[uniform avoidance coupling]
A \emph{uniform avoidance coupling} of two Markov chains $X$ and $Y$ on a common state space $M$ is a Markovian avoidance process of $X$ and $Y$ on $M\times M$ such that for any state $(X_t,Y_t)\in M\times M$ with positive probability in the stationary distribution, the probability of each token transitioning to each other state is the same as the corresponding transition probability in the original chain: $$\sum_{Y_{t+1}\in M} T[(X_t,Y_t)\rightarrow (X_{t+1},Y_{t+1})]=T_X[X_t\rightarrow X_{t+1}],$$ $$\sum_{X_{t+1}\in M} T[(X_t,Y_t)\rightarrow (X_{t+1},Y_{t+1})]=T_Y[Y_t\rightarrow Y_{t+1}].$$
\end{definition}

\begin{theorem}
Any uniform avoidance coupling is an avoidance coupling.
\end{theorem}

\noindent\emph{Proof:} The chains $X$ and $Y$ as described above are each faithful to their original probability matrix, therefore this is a coupling of $X$ and $Y$. \hfill$\Box$\bigskip

\noindent This is a stronger condition than the avoidance coupling simply being Markovian. We are now requiring that not only should there be no correlation between a token's history and its next step, but also the current position of the other token, provided it's a state with positive probability.

\begin{definition}[uniform avoidance coupling of simple random walks]
A \emph{uniform avoidance coupling of two simple random walks} on a graph $G=(V,E)$ is a Markovian avoidance process on $V\times V$ of simple random walks, such that at each state $(x,y)$ with a positive stationary probability, the probability that token $X$ goes to each neighbor in the next move is $\frac{1}{d(x)}$, and the total probability that $Y$ goes to each neighbor in the next move is $\frac{1}{d(y)}$. The moves of the two tokens are, however, not independent.
\end{definition}

\begin{example}(Example 1.4) On a cycle $C_n$ with $n\geq 4$, put two tokens at any pair of non-adjacent vertices. At each state, with probability 1/2 move them both clockwise, and with probability 1/2 move them both anticlockwise.\end{example}

\noindent It is clear in this example, that the overall probability at each step of a particular token going to each adjacent vertex is the same. Therefore, it is a uniform avoidance coupling of simple random walks. However, the moves of the two tokens are not independent, or otherwise the two tokens would collide or cross in time polynomial in $n$. In fact, if we know the move of one token, the other token's move is fixed. Since the clockwise distance on the cycle remains the same in this process, we will call it \emph{fixed distance coupling.}

Notice that no avoidance coupling of simple random walks on a complete graph can be uniform, since that would require a token being able to step on any vertex currently occupied by another token. In fact, at no point can a uniform avoidance coupling of simple random walks be in a state $(X_t,Y_t)$ where $X_t$ and $Y_t$ are adjacent.

\begin{example}(Avoidance coupling that is not a uniform avoidance coupling)

\noindent Consider the coupling from Example 3.3 on $K_4$. Label vertices 1, 2, 3, 4 and let 1, 4 and 2, 3 form the two clusters. Start two tokens in position $(1, 2)$.

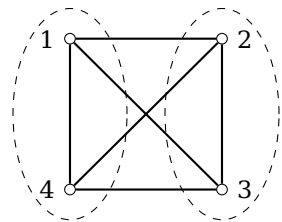
\begin{figure}[ht!]\centering
\begin{tikzpicture}
\draw[thick] (0,0) -- (0,2) -- (2,2) -- (2,0) -- (0,0);
\draw[thick] (0,0) -- (2,2);
\draw[thick] (2,0) -- (0,2);
\foreach \x in {(0,0),(0,2),(2,2),(2,0)}{
\filldraw[white] \x circle (2pt);
\draw \x circle (2pt);
}
\draw[dashed] (0,1) ellipse (0.75cm and 1.4cm);
\draw[dashed] (2,1) ellipse (0.75cm and 1.4cm);
\draw (-0.3,2) node {1};
\draw (2.3,2) node {2};
\draw (2.3,0) node {3};
\draw (-0.3,0) node {4};
\end{tikzpicture}
\end{figure}
\end{example}

\noindent A good way to think about this algorithm is to we consider the moves of the two tokens together. The second token decides which of the available three vertices to move to, and the first token moves out of the way so that at the beginning of the next round it's still in the separate cluster.

Since this process has no memory, at any point the first token is at 1 we don't know if the second token is at 2 or 3, and over all the first token has equal probability of moving to any other vertex. However, if we know the second token is at 2, the first token has probability 2/3 of going to 3 and 1/3 of going to 4. Therefore, this is an avoidance coupling of simple random walks, but not a uniform one.

\begin{remark}
There can be no uniform avoidance coupling of simple random walks on a complete graph. It would require the token that moves first to be able to step on the other token.
\end{remark}

\begin{theorem}
There can be no uniform avoidance coupling of two simple random walks on a tree.
\end{theorem}

\noindent\emph{Proof:} Suppose for contradiction that the coupling exists. Take any state with a non-zero stationary distribution. There is a non-zero chance that $X$ will step towards $Y$, and consequently, there is a non-zero chance that it will step towards $Y$ enough times in a row that $Y$ will have no state to escape to.\hfill$\Box$

\section{Examples of uniform avoidance coupling of simple random walks}

\begin{example}
Let $G$ be the $n$-cube, $n\geq 2$ and start the two tokens in states $(0,\dots,0)$ and $(1,\dots,1)$ respectively. At each step, choose one coordinate uniformly at random and flip it for both tokens.
\end{example}

\begin{example}
Take a bipartite graph where the degree of each vertex is at least 2. Start the chain in a state $(x_0,y_0)$ where $x_0,\ y_0$ are on the same side of the bipartition. Let $x\sim x'$ mean that $x$ and $x'$ are adjacent. Define transition probabilities as follows:\medskip


\noindent Suppose that the chain is in state $(x,y)$ on the same side of the bipartition.\begin{itemize}
\item If $x$ and $y$ have no common neighbors, the transition probability to each state $(x',y')$ such that $x\sim x'$ and $y\sim y'$ is $\frac{1}{d(x)d(y)}$.
\item If $x$ and $y$ have $c$ common neighbors where $c\geq2,$ with the common neighborhood of $x$ and $y$ being $C$, then define transition probabilities as follows:\begin{itemize}
\item If $x',y'\in C,\ x'\not=y'$, then the next state is $(x',y')$ with probability $\frac{c}{d(x)d(y)(c-1)}$.
\item In any other case of $x'\sim x,\ y'\sim ,\ x'\not=y'$, the next state is $(x',y')$ with probability $\frac{1}{d(x)d(y)}$.
\end{itemize}
\item If $x$ and $y$ have exactly one common neighbor $z$, define the transition probabilities as follows. For any $x'\sim x,\ y'\sim y$ where $x',y'\not=z$:\begin{itemize}
\item The next state is $(z,y')$ with probability $\frac{1}{d(x)(d(y)-1)}$.
\item The next state is $(x',z)$ with probability $\frac{1}{d(y)(d(x)-1)}$.
\item The next state is $(x',y')$ with probability $\frac{1-\frac{1}{d(x)}-\frac{1}{d(y)}}{(d(x)-1)(d(y)-1)}=\frac{d(x)d(y)-d(y)-d(x)}{d(x)d(y)(d(x)-1)(d(y)-1)}$

$=\frac{1}{d(x)d(y)}-\frac{1}{d(x)d(y)(d(x)-1)(d(y)-1)}$.
\end{itemize}
\end{itemize}

\noindent We can verify that this is a uniform avoidance coupling of simple random walks as follows. Clearly, for states $(x,y)$ where $x$ and $y$ have no common neighbors, the walks are not only uniform but also independent.

\noindent In the other cases, we want to make sure that for each $x'\sim x$ (an analogous process works for the neighbors of $y$): $$\sum_{y'\sim y}T[(x,y)\rightarrow(x',y')]=\frac{1}{d(x)}$$ And that for each $y'\sim y:$ $$\sum_{x'\sim x}T[(x,y)\rightarrow(x',y')]=\frac{1}{d(y)}$$ If $x$ and $y$ have  a unique, common neighbor $z$, then: $$\sum_{y'\sim y}T[(x,y)\rightarrow(z,y')]=\frac{1}{d(x)(d(y)-1)}\times (d(y)-1)+0=\frac{1}{d(x)}$$ 

\noindent And for any $x'\sim x,\ x'\not=z:$ $$\sum_{y'\sim y}T[(x,y)\rightarrow(z,y')]=\frac{1}{d(y)(d(x)-1)}+\left[\frac{d(x)d(y)-d(y)-d(x)}{d(x)d(y)((x)-1)(d(y)-1)}\right]\times (d(y)-1) $$ $$=\frac{1}{d(y)(d(x)-1)}+\left[\frac{d(x)d(y)-d(y)-d(x)}{d(x)d(y)(d(x)-1)}\right] =\frac{d(x)+d(x)d(y)-d(y)-d(x)}{d(x)d(y)(d(x)-1)}$$ $$=\frac{d(y)(d(x)-1)}{d(x)d(y)(d(x)-1)}=\frac{1}{d(x)} $$

\noindent If the size of the common neighborhood $C$ of $x$ and $y$ is $c\geq 2$, then we can verify the following. For any $x'\sim x$, such that $x'\in C$: $$\sum_{y'\sim y}T[(x,y)\rightarrow(z,y')]=\frac{c}{d(x)d(y)(c-1)} \times(c-1)+\frac{1}{d(x)d(y)}\times(d(y)-c)=\frac{1}{d(x)}$$ And for $x'\sim x$ such that $x'\notin C$, we get: $$\sum_{y'\sim y}T[(x,y)\rightarrow(z,y')]=\frac{1}{d(x)d(y)}d(y)=\frac{1}{d(y)}.$$
This construction is important, because unlike the previous examples,  it is a uniform avoidance coupling that has a higher entropy than a single simple random walk. For a given move of token $X$, we have many possible moves of token $Y$ and vice versa.
\end{example}

We have shown that uniform avoidance coupling of simple random walks is impossible on trees, and that we can construct such coupling for a cycle $C_n,\ n\geq 4$, an $n$-cube with $n\geq 2$ and any bipartite graph in which no vertex has degree 1. In the next section, we will show that for cycles, this is the only Markovian avoidance coupling of simple random walks. We will then identify some classes of graphs for which we can construct such couplings.

\section{Fixed distance is the only Markovian avoidance coupling of SRWs on a cycle}

\noindent We will call the uniform avoidance coupling of two simple random walks on a cycle in examples 1.4 and 4.4 a \emph{fixed distance} avoidance coupling, for obvious reasons.

\begin{theorem}[Toy example: $C_4$]The only Markovian avoidance coupling of simple random walks on $C_4$ is a strong avoidance coupling with fixed distance.\end{theorem} 

\noindent \emph{Proof:} As soon as $X$ and $Y$ are not neighbors at the beginning of a move, the movement of $X$ forces the movement of $Y$. Any avoidance coupling where they remain neighbors forever is not a coupling of simple random walks, since they would go only in one (clockwise or anticlockwise) direction.\hfill$\Box$

\begin{theorem}[Toy example: $C_5$] The only Markovian avoidance coupling of simple random walks on $C_5$ is a uniform avoidance coupling with fixed distance.\end{theorem}

\begin{wrapfigure}[6]{r}[-10pt]{3.7cm}\centering
\raisebox{0pt}[\dimexpr\height-1.3\baselineskip\relax]{\begin{tikzpicture}[scale=0.85]
\foreach \x in {(0,1.6),(-1.5,0.3),(1.5,0.3),(-0.9,-1.3),(0.9,-1.3)}
	\filldraw \x circle (2pt);
\draw (-1.5,0.3) --(0,1.6) --  (1.5,0.3) -- (0.9,-1.3) -- (-0.9,-1.3) -- (-1.5,0.3);
\draw (0,1.8) node {1};
\draw (-1.7,0.3) node {5};
\draw (1.7,0.3) node {2};
\draw (1,-1.55) node {3};
\draw (-1,-1.5) node {4};
\end{tikzpicture}}
\end{wrapfigure}

\noindent \emph{Proof:} First, if the clockwise distance is even, it will stay even and vice versa. Suppose then that we are in the odd case (otherwise consider anticlockwise distance instead). Let us number the vertices $\{1,2,3,4,5\}$ in clockwise order, and let, again, $s(x,y)$ be the stationary distribution of state $(x,y)$ and $T[(x,y)\rightarrow(x',y')]$ the transition probability that while at $(x,y)$, we move to $(x',y')$.

\medskip

\noindent Since $X$ must be at $1$ exactly 1/5 of the moves, we have $s(1,2)+s(1,4)=1/5$. Half the time, it will then step to the right, but that is only possible if we are in fact at $(1,4)$. So: $$s(1,4)\left[T[(1,4)\rightarrow (2,3)]+T[(1,4)\rightarrow(2,5)]\right]=\frac{1}{10}.$$ Now, consider all instances when $X$ goes form 5 to 1 to 2. This happens 1/20 of the time, and is only possible if $Y$ starts at 3, and then goes to 4. So: $$s(5,3)T[(5,3)\rightarrow(1,4)]\left[T[(1,4)\rightarrow (2,3)]+T[(1,4)\rightarrow(2,5)]\right]=\frac{1}{20},$$ $$s(5,3)T[(5,3)\rightarrow(1,4)]=\frac{s(1,4)}{2}.$$ 

\noindent In other words, half the time that the tokens are in state $(1,4)$, they have come there from $(5,3)$. Similarly, we can show that the other half the time they come from $(2,5)$, so the only moves to $(1,4)$ are ones that preserve distance. The same can be shown for all other states at clockwise distance 3, which concludes the proof.\hfill $\Box$\bigskip

\begin{theorem}
Any uniform avoidance coupling on a cycle has fixed clockwise distance from $X$ to $Y$.
\end{theorem}\medskip

\noindent\emph{Proof:} Suppose that the two tokens are at a state with minimum clockwise distance from $X$ to $Y$ of any state with non-zero stationary distribution. Then, there must be 1/2 chance that $X$ steps clockwise and so $Y$ must step clockwise as well or the distance is decreased. Also, there is a 1/2 chance that $Y$ steps anticlockwise and then $X$ must step anticlockwise as well. So the distance cannot increase.\hfill $\Box$\bigskip

\begin{lemma}
Suppose we have a Markovian avoidance coupling of two simple random walks on $C_n$, and $(1,y)$ is a position with minimum clockwise distance from $X$ to $Y$, of all such positions with positive stationary distribution. Then the probability that the distance increases in the next move is 0.
\end{lemma}

\noindent \emph{Proof:} Suppose for contradiction, that the transition probability $(1,y)\rightarrow(n,y+1)$, that is the probability that the distance increases, is $\epsilon>0$. Since this is a minimum distance, transition probability $(1,y)\rightarrow(1,y-1)$ is 0. Then either the probability that $X$ steps to the left, i.e. combined transition probability $(1,y)\rightarrow(n,y+1)$ and $(1,y)\rightarrow(n,y-1)$ is at least $\frac{1}{2}+\frac{\epsilon}{2}$, or the analogous probability that $Y$ steps to the right is. Suppose without loss of generality, that it's the probability that $X$ steps to the left. Suppose that based on the history of the movement of $X$, we can infer with at least $1-\frac{\epsilon}{4}$ probability that the Markov chain is in the state $(1,y)$. Then $X$ steps to the left with probability at least: $$(\frac{1}{2}+\frac{\epsilon}{2})(1-\frac{\epsilon}{4})=\frac{1}{2}+\frac{3\epsilon}{8}-\frac{\epsilon^2}{8}>\frac{1}{2},$$ and we would arrive at a contradiction. 

We will show that for any $\epsilon$, there exists a natural number $K$ such that if $X$ makes $K$ consecutive clockwise circles starting and ending at position 1, then at the end of this sequence $Y$ is at a minimum clockwise distance from $1$ with probability at least $1-\frac{\epsilon}{4}$.


If the coupling is at distance $l$ that is not a minimum distance, then there exists a position $(x_l,x_l+l)$ such that the ratio of transition probability $(x_l,x_l+l)\rightarrow(x_l+1,x_l+l-1)$ to the sum of transition probabilities $(x_l,x_l+l)\rightarrow(x_l+1,x_l+l-1)$ and $(x_l,x_l+l)\rightarrow(x_l+1,x_l+l+1)$ is $\alpha_l>0.$ Let $\alpha$ be the smallest of all $\alpha_l$ such that non-minimum distance $l$ has positive stationary distribution. This is a lower bound on the probability that the clockwise distance decreases as $X$ goes clockwise around the cycle. Then as token $X$ goes around the cycle starting at 1, $Y$ is either at minimum distance or the distance decreases with probability at least $\alpha$. Let $K$ be such that: $$(1-\alpha)^K<\frac{\epsilon}{2n}.$$ Then, as $X$ makes $K$ revolutions starting at 1, regardless of the starting position of $Y$, the final position is $(1,y)$ with probability at least $1-\frac{\epsilon}{4}.$
\hfill$\Box$

\begin{theorem}
Any Markovian avoidance coupling of simple random walks on a cycle $C_n,\ n\geq 4$ is a uniform avoidance coupling with fixed distance.
\end{theorem}

\noindent\emph{Proof:} Lemma 6.4 applies to any position $(x,y)$ at minimum distance, so once the tokens reach a position of minimum distance, the distance can never increase.\hfill$\Box$\bigskip

\noindent It is worth noticing that the initial finite time segment of the process could look differently. But once the process hits a state at a minimum distance, it will never increase the distance again. However, as long as tokens only transition to adjacent vertices, the behavior over any finite time segment will not tell us whether the walks behave like simple random walks.

\section{A network flow condition equivalent to graphs admitting UAC}

\noindent We will sometimes abbreviate uniform avoidance coupling as \emph{UAC}. Unless otherwise stated, we refer to uniform avoidance coupling of two simple random walks.

\begin{definition}[Forbidden States]
A \emph{forbidden state} for a graph $G=(V,E)$ is a state $(x,y)\in V\times V$ such that it is impossible to set up a uniform avoidance coupling of two simple random walks on $G$ in which the stationary probability associated to $(x,y)$ is non-zero. Let $F$ be the set of forbidden states.
\end{definition}

\noindent All states $(x,y)$ such that $x\sim y$ are forbidden, since otherwise the token that moves first would need to have a no-zero probability of stepping on the other token. The fact that we eliminate states $(x,y)$ such that $x\sim y$ is very important. Now, it no longer matters that the moves of the tokens are alternating. The avoidance requirement becomes symmetric with respect to the two tokens.

\begin{example}
On a cycle $C_n$, a state $(x,y)$ is forbidden if and only if $x\sim y$. For any other state, we have a fixed distance coupling that includes these states.
\end{example}

\noindent Call the set of pairs $(x,y)$ such that $x\sim y$ the \emph{zeroth generation forbidden states}, and call this set $F_0$. Then, consider the following recursive relation.

\noindent Having set $F_i$, generate the set $F_{i+1}$ as follows:\begin{itemize}
\item $F_i\subseteq F_{i+1}$
\item For each $(x,y)\notin F_i$ do the following test. Suppose without loss of generality that $d(y)\geq d(x)$. Set up a directed bipartite graph, with $d(x)$ vertices on the left corresponding to the neighborhood $N(x)$, and $d(y)$ vertices on the right corresponding to the neighborhood $N(y)$. For any vertex that is adjacent to both $x$ and $y$, make a copy of it on both sides. Draw an arc from any $x'\in N(x)$ to $y'\in N(y)$ if and only if $(x',y')\notin F_i.$  Do not draw an arc between the two copies of the same vertex.

The state $(x,y)$ passes the test if and only if for any subset $A\subseteq N(x)$ the neighborhood of that subset has size at least $|A|\frac{d(y)}{d(x)}$, that is, by Hall's Theorem, if and only if it is possible to set up flow in this network such that each node on the left produces $1/d(x)$, and each node on the right receives $1/d(y).$ If $(x,y)$ does not pass, add it to $F_{i+1}$.\end{itemize} If the state passes the test, we can allocate transition probabilities into neighboring states in a way that satisfies the definition of uniform avoidance coupling - each pair gets the probability equal to the flow on the arc between them.

\begin{example} Suppose that vertex $x$ has degree 1, and $x\sim y$. Then for any $z$ such that $y\sim z$, $(x,z),(z,x)\in F_1$.\end{example}

\begin{theorem} If $(x,y)$ belongs to the closure $F=\bigcup_{i=0}^{\infty} F_i$, then it is a forbidden state.\end{theorem}

\noindent \emph{Proof:} Clearly, all states in $F_0$ are forbidden states. Suppose for contradiction that a uniform avoidance coupling on $G$ exists, and $i$ is the smallest integer such that a state $(x,y)\in F_i$ exists that has a non-zero stationary distribution. Then, suppose the Markov chain is in this state. The neighboring states that made $(x,y)$ fail the test are forbidden states. So there is no way to distribute the transition probabilities to non-forbidden states to preserve the uniform avoidance coupling.\hfill $\Box$

\begin{example} In a tree, all states are forbidden.\end{example}

\begin{theorem}
If $F\not=V\times V$, a uniform avoidance coupling can be defined on $G$.
\end{theorem}

\noindent \emph{Proof:} Suppose that $(x,y)\notin F$. Then define the Markov chain as in the test above and start it in $(x,y)$. It is a uniform avoidance coupling. \hfill $\Box$\medskip

\noindent Therefore, $F\not=V\times V$ is a necessary and sufficient condition for a uniform avoidance coupling to be possible on $G$.

\begin{theorem}
The forbidden state analysis is a polynomial-time algorithm in $n$ and $m$.
\end{theorem}

\noindent \emph{Proof:} The algorithm for testing whether a graph admits a uniform avoidance coupling runs as follows:\begin{enumerate}
\item Let $\mathcal{F}$ be all pairs of vertices in $V\times V$ that are adjacent, or that are a vertex paired with itself.
\item Let $\mathcal{F'}=\mathcal{F}$.
\item For all pairs $(x,y)$ in $V\times V$ that are not in $\mathcal{F'}$: construct the following network. Take $x$, the neighborhood $N(x)$ of $x$, the neighborhood $N(y)$ of $y$ (if a vertex appears in both, put a ghost copy of that vertex in both sets), and $y$. $x$ will be the source and $y$ the sink. Put an arc of capacity $|N(y)|$ between any pair of vertices in $N(x)$ and $N(y)$ that are currently in $\mathcal{F'}$. Let $l=lcm(|N(x)|,|N(y)|)$. Add $\frac{l}{|N(x)|}$ arcs of capacity 1 from $x$ to any vertex in $N(x)$. Add $\frac{l}{|N(y)|}$ arcs of capacity 1 from any vertex in $N(y)$ to $y$. Find max flow, which will take polynomial time. Unless the max flow is $l,$ add $(x,y)$ to $\mathcal{F'}$.
\item \begin{itemize}\item If $\mathcal{F'}=\mathcal{F}\not= V\times V$, conclude that graph admits UAC. \item Else, if $\mathcal{F'}= V\times V$, conclude graph does not admit UAC. \item Else, let $\mathcal{F}:=\mathcal{F'}$ and go back to (b).\end{itemize}
\end{enumerate}

\noindent The number of vertex pairs to check is initially ${n \choose 2}-m$, where $n$ is the number of vertices and $m$ is the number of edges. It has to decrease at every stage, and so the number of times we need to perform the max flow check cannot exceed: $$\sum_{i=0}^{{n\choose 2}-m}{n\choose 2}-m-i=\frac{1}{2}({n\choose 2}-m)({n\choose 2}-m+1)=O(n^2).$$ The integer max flow algorithm runs in time $O(ml)\simeq O(m^2)$. And so we conclude that this algorithm runs in polynomial time in $n$ and $m$.\hfill$\Box$\bigskip

\begin{remark} As soon as for some vertex $x$ all pairs involving that vertex are eliminated, we know that the graph will fail the forbidden states test. No pairs involving vertices adjacent to $x$ can pass the next round, and no vertices adjacent to those can pass the one after that. This suggests that the above bound on the runtime can be improved.\end{remark}

\begin{remark} Any graph that has a vertex of degree $n{-}1$ will fail the forbidden states test.\end{remark}

\section{Classes of graphs that admit uniform avoidance coupling}

\subsection{Graphs with automorphisms and Cayley graphs}

\begin{theorem} Any graph $G$ that has an automorphism $\varphi: V\rightarrow V$ such that for any $v\in V$, $\varphi(v)\not=v$ and $\varphi(v)\not\sim v$ admits a uniform avoidance coupling.
\end{theorem}\medskip

\noindent \emph{Proof:} Start in a state $(x,y)$ such that $y=\varphi(x)$. Then let $X$ perform a simple random walk on $V$, and whenever $X$ goes to $x'$, let $Y$ go to $\varphi(x')$.\hfill$\Box$\medskip

\noindent  Clearly, such a graph also admits a super-Markovian avoidance coupling. In such case, $X$ performs a simple random walk, and when it's $Y$'s turn, and $X$ is at some position $x$, $Y$ goes to $\varphi(x)$.\medskip

\begin{example} Let $G$ be the $n$-cube, $n\geq 2$ and label the vertices with binary strings in the usual way. Start two tokens in states $00\dots0$ and $11\dots1$ respectively. At each step, choose one coordinate uniformly at random and flip it for both tokens.
\end{example}\medskip

\begin{example}There exist graphs where some vertices have degree 1, and UAC is possible, such as this one:\end{example}
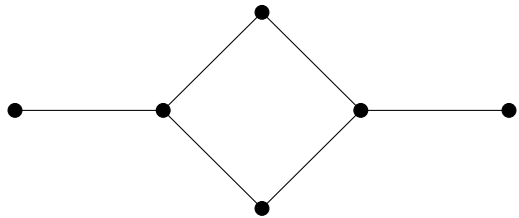
\begin{figure}[ht!]\centering
\begin{tikzpicture}[scale=1.3]
\foreach \x in {(-0.5,0),(1,0),(2,1),(2,-1),(3,0),(4.5,0)}
	\filldraw \x circle (2pt);
\draw (-0.5,0) -- (1,0) -- (2,1) -- (3,0) -- (4.5,0);
\draw (1,0) -- (2,-1) -- (3,0);
\end{tikzpicture}
\end{figure}

\begin{example}
A graph consisting of vertices $\{1,2,\dots,n,1',2',\dots,n'\}$ where the vertices $\{1,2,\dots,n\}$ and $\{1',2',\dots,n'\}$ each form a clique $K_n$, and for each $i\in [n]$, $i\sim i'$. Define the automorphism as $\varphi(i)=(i+1)'$, $\varphi(i')=i+1$. This shows that there exist graphs allowing for UAC with subgraphs $K_n$.
\end{example}\medskip

\begin{remark} Many Cayley graphs satisfy the requirements of Theorem 8.1. Take for example a Cayley graph of the symmetric group $S_n$, $n\geq 3$ with self-inverse generators $\sigma_1,\dots, \sigma_{n-1}$, with $\sigma_i$ corresponding to the transposition of the $i$th  and $i+1$st element in a permutation. Then: $$\varphi(g)=\sigma_1\sigma_2 g $$ defines a uniform avoidance coupling.\end{remark}

\subsection{Regular and strongly regular graphs that admit UAC}

Consider a regular graph of degree $n-2$ with $n$ vertices. In such a graph, for every vertex $x$ there is exactly one vertex $y$ that $x$ is not adjacent to. On such a graph, we can find the following uniform avoidance coupling. Start the process in any state $(x,y)$ such that $x\not\sim y$. At each state, move the token $X$ to any adjacent vertex with the same probability. Move token $Y$ to the vertex that is not adjacent to the vertex $X$ moved to.

\begin{figure}[ht!]\centering
\begin{tikzpicture}
\foreach \x in {(0,1.7),(0,-1.7),(1.7,1),(-1.7,1),(1.7,-1),(-1.7,-1)}{
\foreach \y in {(0,1.7),(0,-1.7),(1.7,1),(-1.7,1),(1.7,-1),(-1.7,-1)}{
\draw \x -- \y;
}
}
\draw[red, thick] (0,1.7) -- (0,-1.7);
\draw[red, thick] (1.7,1) -- (-1.7,-1);
\draw[red, thick] (-1.7,1) -- (1.7,-1);
\foreach \x in {(0,1.7),(0,-1.7),(1.7,1),(-1.7,1),(1.7,-1),(-1.7,-1)}{
\filldraw \x circle (2pt);
}
\end{tikzpicture}
\caption{Edges in the complement are marked in red.}
\end{figure}
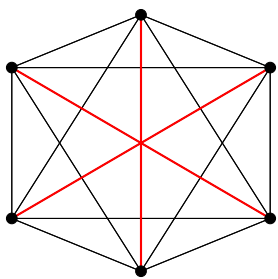

\noindent A similar process can be constructed for any regular graph of degree $n-3$, with $n$ vertices. The compliment of such graph is a collection of cycles. Choose a direction on each of these cycles that will be labelled as \textquotedblleft clockwise." Then start the tokens in a state $(x,y)$ such that $y$ is the clockwise neighbor of $x$ on $x$'s cycle in the complement. Then, have token $X$ perform a simple random walk, and token $Y$ always go to the vertex that in the clockwise neighbor of $X$'s position in the complement.

\begin{remark} Both of the processes described above are also super-Markovian. They are minimum entropy couplings, since one of the tokens is performing a simple random walk and the other's moves are completely determined by what the first one does. \end{remark}

\noindent Not all non-complete regular graphs admit a uniform avoidance coupling.

\begin{example} The graph in Fig. 5 is a 4-regular graph on 12 vertices that does not admit a uniform avoidance coupling.\end{example}

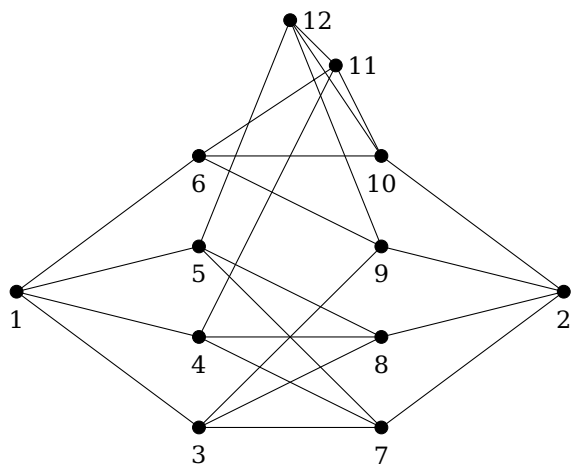
\begin{figure}[ht!]\centering
\begin{tikzpicture}[scale=1.2]
\foreach \x in {0,1,2}{
\filldraw (0,\x) circle (2pt);
\draw (0,\x) -- (-2,1.5);
\foreach \y in {0,1}{
\draw (0,\x) -- (2,\y);
}
}
\foreach \x in {0,1,2}{
\filldraw (2,\x) circle (2pt);
\draw (2,\x) -- (4,1.5);}
\draw (-2,1.5) -- (0,3);
\draw (4,1.5) -- (2,3);
\filldraw (0,3) circle (2pt);
\filldraw (2,3) circle (2pt);
\filldraw (-2,1.5) circle (2pt);
\filldraw (4,1.5) circle (2pt);
\filldraw (1,4.5) circle (2pt);
\filldraw (1.5,4) circle (2pt);
\draw (-2,1.2) node {1};
\draw (4,1.2) node {2};
\draw (0,-0.3) node {3};
\draw (0,1.7) node {5};
\draw (0,2.7) node {6};
\draw (0,0.7) node {4};
\draw (2,-0.3) node {7};
\draw (2,0.7) node {8};
\draw (2,1.7) node {9};
\draw (2,2.7) node {10};
\draw (1.8,4) node {11};
\draw (1.3,4.5) node {12};
\draw (0,0) -- (2,2);
\draw (0,1) -- (1.5,4);
\draw (0,2) -- (1,4.5);
\draw (2,3) -- (1,4.5);
\draw (2,3) -- (1.5,4);
\draw (0,3) -- (2,3);
\draw (0,3) -- (2,2) -- (1,4.5) -- (1.5,4) --(0,3);
\end{tikzpicture}
\label{graph12}
\caption{Clearly $\mathcal{F}\not=F_0$, since the non-adjacent pair $(1,2)\in F_1$. A forbidden state analysis on this graph eliminates all pairs.}
\end{figure}

\begin{lemma}
A regular graph admits a uniform avoidance coupling of simple random walks if and only if it admits a minimum entropy one.
\end{lemma}

\noindent\emph{Proof:} Suppose that a regular graph $G$ passes the forbidden states test. Then for any pair of vertices $(x,y)\not\in\mathcal{F}$ there exists a matching between their neighborhoods such that a vertex $v\in N(x)$ is matched with some vertex $w$ where $(v,w)\not\in\mathcal{F}$. Pick such matching for each pair $(x,y)\not\in\mathcal{F}$. 

These define a minimum entropy UAC of SRWs, in which one token performs a simple random walk and token $Y$ steps to the state defined by the appropriate matching.\hfill$\Box$

\begin{definition}
A graph $G=(V,E)$ with $n$ vertices is \emph{strongly regular} if there exist natural numbers $k,\ \lambda,\ \mu$ such that:\begin{itemize}\item the degree of every vertex is $k$, \item any two adjacent vertices have exactly $\lambda$ common vertices, \item any two non-adjacent vertices have exactly $\mu$ common neighbors.
\end{itemize}
\end{definition}

Suppose that $G=(V,E)$ is a connected, non-complete strongly regular graph with parameters $(n,k,\lambda,\mu)$. We will show that unless $\mu>k/2$, not only does it admit a uniform avoidance coupling of simple random walks, but a stronger statement is true. In fact, $F_0=F$ i.e. the only forbidden states are $(x,y)$ where $x$ and $y$ are adjacent.

\begin{theorem}If $G$ is a connected, non-complete strongly regular graph with $\max(\lambda,\mu)\leq k/2$, or with $\lambda<k/2$, then on $G$ we have $F=F_0$.\end{theorem}

\begin{corollary}A uniform avoidance coupling is possible on connected, non-complete strongly regular graphs with $\max(\lambda,\mu)\leq k/2$, as well as those with $\lambda<k/2$.\end{corollary}

\begin{remark} This includes Conference Graphs and many other classic examples of strongly regular graphs.\end{remark}

\noindent \emph{Proof of Thm 8.10:} To prove this, we need to show that no non-adjacent pair is added to $F_1$. Since the graph is regular, this is equivalent to showing that for any $x,\ y\in V$ such that $x\not\sim y$, we can find a perfect matching between the neighbourhoods such that vertices are only matched to non-adjacent vertices. Then when the process is at state $(x,y)$, and token $X$ goes to some $x'$, token $Y$ goes to the neighbor $y'$ of $y$ that $x'$ is matched to. Since these are not adjacent, we can continue the process.

\begin{figure}[ht!]\centering
\hfill
\begin{tikzpicture}[scale=0.75]
\foreach \x in {5,4,3}{
\draw[thick,black!30] (0,\x) -- (2.5,\x);
}
\foreach \x in {0,1,2,3,4,5}{
\filldraw (0,\x) circle (2pt);
\filldraw (2.5,\x) circle (2pt);
\draw (0,\x) -- (-2,2.5);
\draw (2.5,\x) -- (4.5,2.5);
\draw (-2,2.1) node {$x$};
\draw (4.5,2.1) node {$y$};
}
\filldraw (-2,2.5) circle (2pt);
\filldraw (4.5,2.5) circle (2pt);
\draw[thick] (0,5) -- (2.5,3);
\draw[thick] (0,5) -- (2.5,0);
\draw[thick] (2.5,5) -- (0,3);
\draw[thick] (0,4) -- (2.5,2);
\draw[thick] (0,4) -- (2.5,1);
\draw[thick] (0,3) -- (2.5,2);
\foreach \x in {4, 1, 0}{
\draw[thick] (0,0) -- (2.5,\x);
}
\foreach \x in {5, 1, 0}{
\draw[thick] (0,1) -- (2.5,\x);
}
\foreach \x in {4, 3, 2}{
\draw[thick] (0,2) -- (2.5,\x);
}
\draw[dashed] (0,4) ellipse (0.75cm and 1.3cm);
\draw[dashed] (2.5,4) ellipse (0.75cm and 1.3cm);
\draw (0,5.5) node {$N(x)$};
\draw (2.5,5.5) node {$N(y)$};
\end{tikzpicture}\hfill
\begin{tikzpicture}[scale=0.75]
\draw (0,5.5) node {$N(x)$};
\draw (2.5,5.5) node {$N(y)$};
\foreach \x in {0,1,2,3,4,5}{
\filldraw (0,\x) circle (2pt);
\filldraw (2.5,\x) circle (2pt);
}
\foreach \x in {5,3,2}{
\draw[thick,purple] (0,0) -- (2.5,\x);
}
\foreach \x in {4,3,2}{
\draw[thick,purple] (0,1) -- (2.5,\x);
}
\foreach \x in {5,1,0}{
\draw[thick,purple] (0,2) -- (2.5,\x);
}
\foreach \x in {4,1,0}{
\draw[thick,purple] (0,3) -- (2.5,\x);
}
\foreach \x in {5, 3, 1}{
\draw[thick,purple] (0,4) -- (2.5,\x);
}
\foreach \x in {4, 1, 2}{
\draw[thick,purple] (0,5) -- (2.5,\x);
}
\end{tikzpicture}\hfill\hfill
\caption{On the left, an example with $k=6,\ \mu=3,\ \lambda=2$ we get 3 vertices that need to be represented on both sides. The two copies are joined by a \textquotedblleft ghost edge". On the right, the complement graph between the two neighborhoods. This is the graph in which we need to find a perfect matching.}
\label{bipartite}
\end{figure}
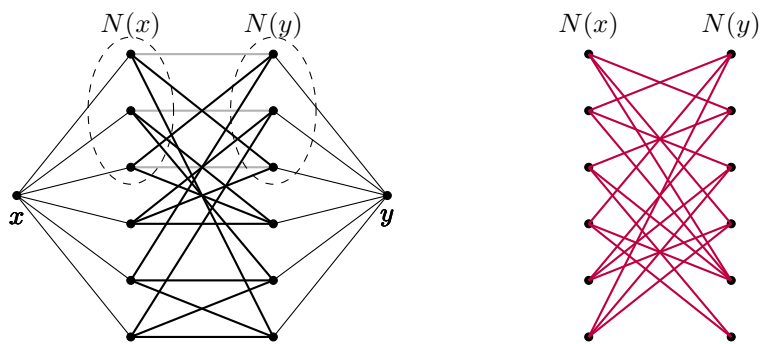

Let us represent the neighborhoods of $x$ and $y$ on two sides. $x$ and $y$ will have exactly $\mu$ neighbors in common. We will put a copy of those on each side and join those copies by a ghost edge (represented in grey in Fig. \ref{bipartite}) to make sure we don't match a vertex to itself. Now, every vertex of those $\mu$ is adjacent to $y$, and so has $\lambda$ common neighbors with $y$. So it will have an additional $\lambda$ edges going to the other side. Every vertex that is not adjacent to $y$ has $\mu$ common neighbors with $y$. The complement of this graph between the two neighborhoods is where we need to find a perfect matching.

In the example of Fig. \ref{bipartite}, we know immediately that we will be able to find a perfect matching, since this bipartite graph is regular. Therefore, if $k=6,\ \mu=3,\ \lambda=2$, or in fact in any case where $\mu=\lambda+1$, we know that $F=F_0$ and so such a graph admits a uniform avoidance coupling. Notice that since the common neighborhood of $x$ and $y$ (marked by the ellipses) are the same vertices appearing on both sides, adjacency relations within that set need to be symmetric.

So, under what circumstances would this argument fail? We would need there to be, for some $x$ and $y$ such that $x\not\sim y$, subset $A$ of the neighborhood $N(x)$ of $x$ and subset $B$ of the neighborhood $N(y)$ of $y$, with $|A|>|B|$, such that every node in $N(y)$ that is not adjacent to at least one vertex in $A$ is in $B$. This is a tall order. Every node in $A$ has to be adjacent to every node in $N(y)-B$, except possibly itself.

Let's take a closer look at what we know about the complement bipartite graph between the two neighborhoods. Now, edges will represent two distinct nodes in the neighborhoods that are not adjacent in the original graph (marked in purple.) There are $\mu$ nodes in common, and each has degree $k-\lambda-1$, since in the original graph we need to exclude the other copy of the same node, and the $\lambda$ nodes in $N(y)$ to which this node is adjacent. Any adjacency relations between the two copies of this set need to be symmetric. The remaining $k-\mu$ nodes have degree $k-\mu$.\medskip

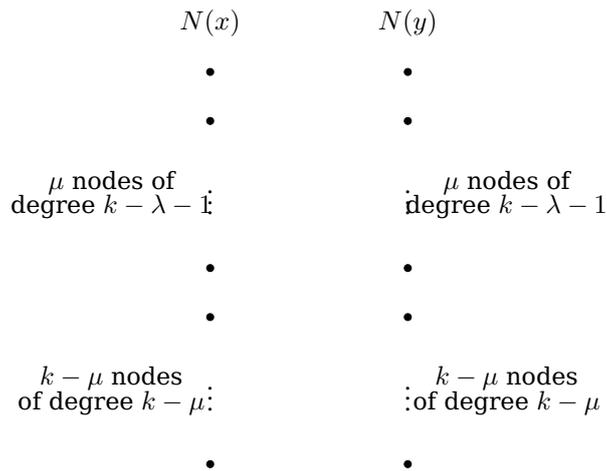
\begin{figure}[ht!]\centering
\begin{tikzpicture}[scale=0.65]
\draw (0,9) node {$N(x)$};
\draw (4,9) node {$N(y)$};
\foreach \x in {0,3,4,7,8}{
\filldraw (0,\x) circle (2pt);
\filldraw (4,\x) circle (2pt);
}
\draw (0,1.5) node {$\vdots$};
\draw (4,1.5) node {$\vdots$};
\draw (0,5.5) node {$\vdots$};
\draw (4,5.5) node {$\vdots$};
\draw (-2,1.75) node {$k-\mu$ nodes}; 
\draw (-2,1.25) node {of degree $k-\mu$};
\draw (6,1.75) node {$k-\mu$ nodes}; 
\draw (6,1.25) node {of degree $k-\mu$};
\draw (-2,5.75) node {$\mu$ nodes of}; 
\draw (-2,5.25) node {degree $k-\lambda-1$};
\draw (6,5.75) node {$\mu$ nodes of}; 
\draw (6,5.25) node {degree $k-\lambda-1$};
\end{tikzpicture}
\caption{Degrees of vertices in the complement graph, in which we are looking for a perfect matching.}
\end{figure}

Then $F_0=F_1$ unless there exist some vertices $x,\ y$, $x\not\sim y$, for which we can find sets $A\subset N(x),\  B\subset N(y)$ with the following properties. Either $|A|>|B|$ and every vertex in $A$ adjacent to every vertex in $N(y)-B$ except possibly itself, or $|A|<|B|$ and every vertex in $B$ is adjacent to every vertex in $N(x)-A$ except possibly itself. Therefore, we can assume without loss of generality that $|A|>k/2$, since otherwise we can take $N(y)-B$ rather than $A$. (Fig. 8)

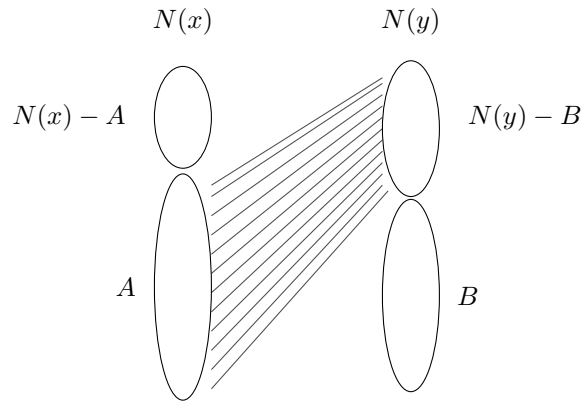
\begin{figure}[ht!]\centering
\begin{tikzpicture}[scale=0.75]
\filldraw[black!70] (0.5,1.8) -- (3.5,3.7);
\draw[black!70] (0.5,-1.8) -- (3.6,1.7);
\foreach \x in {0.2,0.4,0.6,0.8,1,1.2,1.4,1.6,1.8,2}{
\draw[black!70] (0.5,-1.8+1.7*\x) -- (3.5,1.6+\x);
}
\draw (0,0) ellipse (0.5cm and 2cm);
\draw (4,-0.15) ellipse (0.5cm and 1.7cm);
\draw (0,3) ellipse (0.5cm and 0.9cm);
\draw (4,2.8) ellipse (0.5cm and 1.2cm);
\draw (-1,0) node {$A$};
\draw (5,-0.15) node {$B$};
\draw (-2,3) node {$N(x)-A$};
\draw (6,3) node {$N(y)-B$};
\draw (0,4.7) node {$N(x)$};
\draw (4,4.7) node {$N(y)$};
\end{tikzpicture}
\caption{Every vertex in $A$ is adjacent every vertex in $N(y)-B$ except possibly itself.}
\end{figure}

Let $|N(y)-1|>1$, since otherwise $A=N(x)$ and some vertex in $N(y)$ is connected to every vertex in $N(x)$, which would necessarily mean that $\lambda=k$ and the graph would be complete. 

Any two distinct vertices in $N(y)-B$ have at least $|A|$ neighbors in common. This is clear if neither of them is also in $A$, in which case they have at least $|A|+1$ common neighbors, where the 1 comes from $y$. If one of them is also in $A$, this number is at least $|A|$. If they are both in $A$, they have $|A|-2$ common neighbors in $A$, and also both $x$ and $y$, which again makes $|A|$. We conclude that $\max(\lambda,\mu)> k/2$.

Notice that each vertex that is in both $N(x)$ and $N(y)$ has exactly $\lambda$ neighbors in $N(x)$, and every vertex in $N(y)$ that is not in $N(x)$ has exactly $\mu$ neighbors in $N(x)$. Suppose that $\lambda< k/2.$ Then none of the vertices in $N(y)-B$ can be in $N(x)$. But then there are at most  $k-\mu$ of them, so we must have $|A|>\mu$. But then the vertices on the other side have to have degree $|A|$, and we arrive at a contradiction, and conclude that $\lambda\geq k/2.$\hfill $\Box$ \medskip

\begin{remark} In particular, unless for all  $x,\ y$, $x\not\sim y$ all pairs of vertices in $N(y)-B$ are non-adjacent, $\lambda>k/2$. Unless for all  $x,\ y$, $x\not\sim y$ $N(y)-B\subseteq N(x)$, $\mu\geq |A|>k/2$, since then any vertex that is not in $N(x)$ would have at least $|A|$ common neighbors with $x$. Unless all pairs of vertices in $N(y)-B$ are adjacent, $\mu>k/2$.\end{remark}

\begin{remark}
The coupling that this defines is minimum-entropy, but not necessarily super-Markovian. While the walk of $Y$ is determined completely by the walk of $X$, the matching might change based on which vertex $x$ token $X$ moved to $x'$ from.
\end{remark}

\medskip

\noindent Coupling of Markov chains often comes up in the context of the mixing time of Markov chains. It was pioneered by Aldous \cite{Aldous} and Bubley and Dyer\cite{BD}. Couplings of this kind are uniform couplings - in fact, some literature \cite{VG} appears to define coupling as in Definition 1.1 and then assume it's uniform. We have demonstrated in Example 4.5 that this is not necessarily the case.


\ACKNO{I would like to thank my thesis advisor, Prof. Peter Winkler, as well as Drs Alexander Holroyd and Ralf Banisch for their insight. This work was supported by a Dartmouth College graduate student scholarship as well as NSF grants DMS-0901475 and DMS-1162172.}


\end{document}